\documentclass[11pt]{article}
\usepackage{latexsym,amsmath,amsthm,amssymb}
\usepackage[left=1in,top=1in,right=1in,bottom=1in]{geometry}
\usepackage{setspace}

\def\qed{\hfill\ifhmode\unskip\nobreak\fi\quad\ifmmode\Box\else\hfill$\Box$\fi}
\def\ite#1{\hfill\break${}$\hbox to 50pt {\quad(#1)\hfill}}
\newtheorem{thm}{Theorem}

\newtheorem{lem}[thm]{Lemma}

\newtheorem{claim}[thm]{Claim}

\def\eps{\varepsilon}

\begin{document}

\pagestyle{myheadings}
\markright{{\small{\sc F\"uredi, Ruzsa:  Nearly subadditive sequences, 
}}}

\title{\vspace{-0.5in}Nearly subadditive sequences
}

\author{
{{Zolt\'an F\" uredi}}\thanks{
\footnotesize {Alfr\' ed R\' enyi Institute of Mathematics, Hungary
E-mail:  \texttt{furedi.zoltan@renyi.mta.hu}.
Research supported in part by the Hungarian National Research, Development and Innovation Office NKFIH, K116769, and
 by the Simons Foundation Collaboration Grant 317487.
}}
\and
{{Imre Z. Ruzsa}}\thanks{
\footnotesize {Alfr\' ed R\' enyi Institute of Mathematics, Hungary
E-mail:  \texttt{ruzsa.imre@renyi.mta.hu}.
 Research is supported in part  by
ERC–AdG Grant No.321104 and Hungarian National Foundation for Scientific Research Grant 
  NK104183. 
}}
}

\date{\today}

\maketitle

\vspace{-0.3in}

\begin{abstract} 
We show that the de Bruijn-Erd\H{o}s condition for the error term in their
 improvement of Fekete's Lemma is not only sufficient but also necessary
 in the following  strong sense.
Suppose that given a sequence $0\leq f(1)\leq f(2)\leq f(3)\leq \dots $ such that
\begin{equation}\label{eq01}\sum_{ n=1}^{\infty}   f(n)/n^2 = \infty.
   \end{equation}
Then, there exists a sequence  $\{b(n)\}_{n=1,2,\dots}$ satisfying
\begin{equation}\label{eq1} b(n+m)  \leq    b(n) + b(m) + f(n+m)
   \end{equation} 
such that the sequence of slopes $\{ b(n)/n\}_{n=1,2,\dots}$ takes every rational number.

When the series~\eqref{eq01}
is bounded
   we improve their result as follows.
If  there exist $N$ and real $\mu >1$ such that~\eqref{eq1}
  holds for all pairs $(n,m)$ with $N\leq  n\leq m \leq \mu  n$, then
 $\lim_n b(n)/n $ exists. 

\end{abstract}


\section{Fekete's lemma on subadditive sequences}\label{S1} 

An infinite sequence of reals  $a(1), a(2), \dots, a(n), \dots$  is called {\em subadditive} if
\begin{equation} \label{eq03}  a(n+m) \leq a(n)+ a(m)
   \end{equation}
holds for all integers $n,m\geq 1$.

Every (reasonable) calculus textbook contains Fekete's~\cite{MFekete} Lemma as a theorem
(or as an exercise, see, e.g., Polya and Szego~\cite{PSz}).
It says  that if the sequence $\{ a(n) \}$ is subadditive,
then the sequence $\{ a(n)/n\}$  has a limit (possible negative infinity).
Moreover, that limit is equal to the infimum,
 \begin{equation}\label{eq04}  \lim_{n\to \infty}\dfrac{a(n)}{n} = \inf_{k\geq 1}
 \frac{a(k)}{k}.
   \end{equation}

\pagebreak

{\bf The standard proof of Fekete's Subadditive Lemma}

Using the subadditivity we get by induction (from $n-k$ to $n$)  that\
\[    a(n) \leq a(k)+a(n-k) \leq 2a(k)+ a(n-2k)\leq \dots  \leq   \lfloor n/k \rfloor a(k) + a(\beta) , 
   \]
 where $0\leq \beta\leq k-1$. (We may define $a(0)=0$).
This implies that  for all $n\geq k\geq 1$
\begin{equation}\label{eq05}
{\frac{a(n)}{n} \leq \frac{a(k)}{k} } +\frac{\max\left\{ |a(1)|, \dots, |a(k-1)|\right\}}{n}.
    \end{equation}
Therefore
\begin{equation}\label{eq06}
\limsup \frac{a(n)}{n} \leq \frac{a(k)}{k}.
   \end{equation}
This holds for every $k$, so
\begin{equation}\label{eq07}
  \limsup \leq \inf,
   \end{equation}
implying\enskip $\limsup=\inf$,  so the limit exists.  \qed

\smallskip
{\bf  A remark on large values of $(n,m)$}

Note that the above proof yields that if the subadditivity~\eqref{eq03} only holds for $n,m \geq N$, then the limit still exists.
We have $a(n)/n \leq a(k)/k$ for all $n\geq k\geq N$ whenever $n/k$ is an integer.
In general, we use induction only if both $k$ and $n-k$ is at least $N$, i.e., we choose $\beta\in [k+1, 2k-1]$.
Instead of~\eqref{eq05} we obtain that for all $n\geq 2k$, $k\geq N$
\begin{equation}\label{eq08}
{\frac{a(n)}{n} \leq \frac{a(k)}{k} } +\frac{\max\left\{ |a(k+1)|, \dots, |a(2k-1)|\right\}}{n}.
    \end{equation}
This implies~\eqref{eq06} and~\eqref{eq07} for $k\geq N$. We obtain
 \begin{equation}\label{eq09}  \lim_{n\to \infty}\dfrac{a(n)}{n} = \inf_{k\geq N}
 \frac{a(k)}{k}.  \quad\qed
   \end{equation}

\smallskip
{\bf  Having the threshold $N$ is a true (and nontrivial) extension}

One might be tempted to think that~\eqref{eq09} can be easily obtained from the original Fekete's lemma~\eqref{eq04}.
Maybe so, but let us consider the following sequence.
Suppose that $2\leq N \leq n_1 < n_2 < n_3< \dots$ are integers such that $n_{i} - N \leq n_{i+1}$.
Define for all $i\geq 1$ and positive integer $n$
$$a(n):=\begin{cases} 1 \qquad \qquad\quad n\le n_1\\
             1 \qquad\qquad\quad n_{i+1}-N  \leq n\leq n_{i+1}-2,\\
    n/n_i \qquad\quad\,\, \, n_{i}\leq n < n_{i+1}\enskip  ({\rm but }\enskip  |n-n_{i+1}|\notin [2,N]).\end{cases}$$
This sequence satisfies subadditivity for $m,n\geq N$.
Suppose that $\limsup n_{i+1}/n_i = \infty$.
Then the sequence $\{ a(n) \}$ does not seem to be easily transformed to a true subadditive one, because
there are infinitely many  $(x,y)$ pairs with $1\leq x< N$ and $x+y= n_{i+1}-1$ such that
  $a(x+y)-a(y)-a(x)= (n_{i+1}-1)/n_i -2$ is arbitrarily large.

(If one prefers an integer sequence, then can observe that $\{ \lceil a(n)\rceil\}$ has the same properties).

\section{Sub-2 sequences by de Bruijn and Erd\H os}\label{S2}

A sequence $\{ a(n) \} $ is called $\mu$-{\em subadditive} with a {\em threshold} $N$ ($(\mu,N)$-subadditive, for short) if
\begin{equation} \label{eq110}  a(n+m) \leq a(n)+ a(m)
   \end{equation}
holds for all integers $n,m$ such that
 \begin{equation} \label{eq111}    N\leq n \leq m \leq \mu n.
   \end{equation}

\begin{thm}[de Bruijn and Erd\H os, Theorem 22. in~\cite{deBE_1952}]\label{th21}
Suppose that the sequence $\{ a(n) \}$ satisfies~\eqref{eq110} for all integers $N\leq n\leq m\leq 2n$.
Then the sequence of slopes $\{ a(n)/n\}$  has a limit (possible negative infinity).
Moreover, that limit is equal to the infimum,
 \begin{equation*} 
   \lim_{n\to \infty}\dfrac{a(n)}{n} = \inf_{k\geq N} \frac{a(k)}{k}.
   \end{equation*}
  \end{thm}

Actually, they considered the case $N=1$ only.
Here we present a greatly simplified proof.

\smallskip
{\bf  A new proof for Theorem~\ref{th21}}

Fix a $k$, $k\geq N$.
Write $n$ as $n= \left(\lfloor n/k \rfloor-1\right) k + \beta$
where $k \leq \beta\leq 2k-1$.
We will show that
 \begin{equation}\label{eq113}
  a(n)\leq \left(\lfloor n/k \rfloor-1\right) a(k) + a(\beta).
   \end{equation}
This implies that for all $n\geq 2k$, $k\geq N$ inequality~\eqref{eq08} holds, implying~\eqref{eq09} as in earlier proofs, and we are done.

To prove~\eqref{eq113} we need a definition.
A sequence of (positive) integers $X:= \{ x_1, x_2, \dots, x_t\}$ (here $t\geq 1$) is called 2-{\em good} if
  $1/2 \leq x_i/x_j\leq 2$ holds for all $1\leq i, j \leq t$.
If $X$ is a 2-good sequence of length $t$ and we take
 two minimal members, $x_i, x_j\in X$, delete them from $X$ but join $x_{new}:= x_i+x_j$, then the new sequence $X':=X\setminus \{ x_i, x_j\} \cup \{ x_{new}\}$ is 2-good as well.
Note that the sum of the members of $X$ is the same as in $X'$.
If the sequence of $\{ a(x)\}$ is 2-subadditive then   $a(x_{new})\leq  a(x_i)+a(x_j)$ implies that
 \begin{equation}\label{eq114}  \sum_{x\in {X'}} a(x) \leq \sum_{x\in X} a(x).
   \end{equation}

Define the set $X_{\lfloor n/k \rfloor}$ of length $\lfloor n/k \rfloor$ as $\{ k,k,k, \dots, k, \beta\}$.
It is obviously a  2-good sequence with sum $n$.
Define the sets $X_t$ of length $t$ for ${\lfloor n/k \rfloor}\geq t\geq 1$ by the above rule, $X_{t-1}:=X_t'$.
We obtain
$X_{\lfloor n/k \rfloor}$  $\longrightarrow\dots X_t \longrightarrow X_{t-1}\longrightarrow \dots \longrightarrow X_1=\{ n\}$. Then~\eqref{eq114} gives
\[    a(n) = \sum_{x\in X_1} a(x) \leq \dots \leq  \sum_{x\in X_{t}} a(x) \leq \dots \leq \sum_{x\in X_{\lfloor n/k \rfloor}} a(x) =  \left(\lfloor n/k \rfloor-1\right) a(k) + a(\beta) .
\quad \qed
   \]

\section{Sub-$\mu$ sequences with $\mu<2$}\label{S3}

Concerning their result (Theorem~\ref{th21} above) de Bruijn and Erd\H os~\cite{deBE_1952} state, maybe somewhat carelessly, that 'It may be remarked that the inequality in $(7.1)$ cannot be replaced by $\mu^{-1}n \leq m \leq \mu n$ for any $\mu < 2$'.
In their papers~\cite{deBE_1951, deBE_1952} they deal with many conditions and sequences, we could not really know what was in their minds, but our first new result is a strengthening of Theorem~\ref{th21}  for all $\mu>1$.
We show that their condition can be weakened such that the limit exists if~\eqref{eq110} holds only for the pairs $(n,m)$ with $n\leq m\leq \mu n$ for some fixed $\mu >1$.

\begin{thm}\label{th32}
Suppose $\mu>1$ and $N\geq 1 $ are given.
If the sequence $\{ a(1)$, $a(2),\dots \}$ is $(\mu, N)$-subadditive, i.e.,
\[     a(n+m)\leq a(n)+a(m)\quad \forall n\leq m\leq \mu n,\enskip n,m\geq N,
   \] 
then the $\lim_{n\to \infty}\dfrac{a(n)}{n}$ exists and is equal to $\inf_{k\geq N} \dfrac{a(k)}{k}$.
(It may be\enskip  $-\infty$).
\end{thm}

For the proof we investigate sequences $\{ a(n)\} $ where the subadditivity holds only for a very few pairs $(n,m)$.

\smallskip
{\bf  Sub-$\mathbf 1^+$ sequences}

Given $N\geq 1$ a sequence $\{ a(n)\} $ is called $(1^+, N)$ {\em subadditive} if the following two inequalities hold for all $n\geq N$.
\begin{eqnarray*} a(2n) &\leq& a(n)+a(n)\\
    a(2n+1) &\leq&       a(n)+a(n+1).
 \end{eqnarray*}
Given a sequence $\{ a(n)\}$  let
$q(n):=\max \left\{ \dfrac{a(n)}{n},\dots, \dfrac{a(2n-1)}{2n-1}, \dfrac{a(2n)}{2n}\right\}$.

\begin{lem}\label{le33}
Suppose that $N\geq 1$ and the sequence $\{ a(n)\} $ is $(1^+, N)$ subadditive.
Then for $n\geq N$ the sequence  $\{ q(n)\}$ is non-increasing,
   $q(n)\geq q(n+1)$.
  \end{lem}

We only have to show that  $q(n)$ is at least as large as $a(2n+1)/(2n+1)$ and $a(2n+2)/(2n+2)$.
The $1^+$ subadditivity implies
\[
q(n)\geq
\begin{cases} \dfrac{a(n+1)}{n+1}\geq \dfrac{a(2n+2)}{2n+2}, \\
\max \left\{ \dfrac{a(n)}{n}, \dfrac{a(n+1)}{n+1}\right\}\geq
            \dfrac{n}{2n+1} \dfrac{a(n)}{n} +\dfrac{n+1}{2n+1} \dfrac{a(n+1)}{n+1}
                          \geq \dfrac{a(2n+1)}{2n+1}.
\quad \qed
\end{cases}
\]

\smallskip
{\bf  Proof of Theorem~\ref{th32}}

Since the case $\mu\geq 2$ is covered by Theorem~\ref{th21}, we may suppose that $1< \mu< 2$. 
Define the positive integer $k$ by
\begin{equation*} 
    (1+\mu)^{k-1}\leq 2^{k+1} < (1+\mu)^k.
   \end{equation*}
Given any $n$ define the sequences $u_0$, $u_1, \dots $, $u_k$ and $v_0$, $v_1, \dots $, $v_k$
as follows.
\begin{equation*} 
 u_0=v_0:= n, \quad u_{i+1} := 2u_i, \quad  v_{i+1} := v_i + \lfloor \mu v_i\rfloor, \quad (i=0,1,\dots, k-1).
      \end{equation*}
We have $u_k= 2^kn$ and $v_k > (1+\mu)^kn -(1+\mu)^k/\mu$.
So there exists an $N_1$ (depending only from $\mu$) such that $2u_k\leq v_k$ holds in the above process for every integer $n\geq N_1$.

Let $N_2:= \max \{ N, 1/(\mu -1) \}$.
Then the sequence $\{ a(n)\}$ is $(1^+, N_2)$ subadditive.
Lemma~\ref{le33} implies that $L=\lim_{n\to \infty} q(n)$ exists.
If $L= -\infty$ then $\lim_{n\to \infty} a(n)/n = -\infty $ as well, and we are done.
Since $L< \infty$, from now on, we may suppose that  $L$ is a real number.

Choose an (arbitrarily small) $\varepsilon>0$.
There exists an $N_3$ (depending on $\varepsilon$, $\mu$, $N$, and $\{ a(n)\}$)  such that
  $q(n)< L+\varepsilon$ for every $n\geq N_3$.
By the definition of $q$ we get
\begin{equation}\label{eq317}
  a(n)/n< L+\varepsilon
\end{equation}
for every $n\geq N_3$.
We are going to show that for $n\geq \max\{ N_1, N_2, N_3\}$
\begin{equation}\label{eq318}
   a(n)/n> L+ \eps -\varepsilon\left( 1+\mu \right)^k.
\end{equation}
Since this holds for every $\varepsilon >0$ the limit $a(n)/n$ exists and is equal to $L$.

To prove~\eqref{eq318} we need the following claim which holds for each $i\in \{0, 1, \dots, k-1\}$.
\begin{claim}\label{cl34}
If $a(w)/w \leq L+\eps- \eta$ for every $w\in [u_i, v_i]$, then  $a(z)/z < L+\eps- \frac{\eta}{1+\mu}$ for every $z\in [u_{i+1}, v_{i+1}]$.
  \end{claim}
Indeed, every $z\in  [u_{i+1}, v_{i+1}]$ can be written in the form $z=x+y$ where
$x\in [u_i, v_i]$, $x\leq y\leq \mu x$.
Apply subadditivity for $(x,y)$ and the upper bound $L+\eps -\eta$ for $a(x)/x$ and the upper bound $L+\eps$ for $a(y)/y$.
We obtain
\begin{eqnarray*}
                   \frac{a(z)}{z}&=& \frac{a(x+y)}{x+y}\leq \frac{a(x)+a(y)}{x+y}\\
       &=& \frac{a(x)}{x}\frac{x}{x+y} +  \frac{a(y)}{y}\frac{y}{x+y}
    <  (L+\eps -\eta)\frac{x}{x+y} +  (L+\eps)\frac{y}{x+y}
   \\ &=&  L+\eps -\eta \frac{x}{x+y} \leq  L+\eps -\eta \frac{1}{1+\mu}. \qed
   \end{eqnarray*}

\noindent
{\em The end of the proof of Theorem~\ref{th32}}. \enskip
Consider any $n$ with $n\geq \max\{ N_1, N_2, N_3\}$.
By~\eqref{eq317} we have $a(n)/n = L+\eps -h$ for some $h>0$.
Consider the intervals $[u_i,v_i]$ for $i=0,1, \dots, k$, where $[u_0,v_0]$ consists of a single element, namely $n$.
Using Claim~\ref{cl34} we get that $a(x) < L+\eps -h/(1+\mu)^i$ for each $x\in [u_i,v_i]$ for $1\leq i\leq k$.
Especially, $a(x)/x < L+\eps -h/(1+\mu)^k$  for each $x\in [u_k,v_k]$.
Since $2u_k\leq v_k$ we obtain $q(u_k)<  L+\eps -h/(1+\mu)^k$.
But  $q(u_k)\geq L$. This implies $h< \eps (1+\mu)^k$.
We obtained that $a(n)/n= L+\eps -h> L+\eps -\eps (1+\mu)^k$ as claimed in~\eqref{eq318}.
This completes the proof of Theorem~\ref{th32}. \qed

\section{Nearly subadditive sequences, \\ an error term by de Bruijn and Erd\H os}\label{S4}

Let   $f(n)$  be a non-negative, non-decreasing sequence. deBruijn and Erd\H{o}s~\cite{deBE_1952} called the sequence $\{ a(n)\}$ subadditive with an {\em error term} $f$ (or {\em nearly $f$-subadditive}, or $f$-{\em subadditive} for short) if
\begin{equation}\label{eq415}
         a(n+m) \leq a(n)+ a(m) + f(n+m)
         \end{equation}
holds for all positive integers $n,m\geq 1$.
The case $f(x)=0$ corresponds to the cases discussed above.

They showed that if the error term  $f$  is small,
\begin{equation}\label{eq416}
         \sum_{n=1}^{\infty} \,\,  f(n)/n^2  \,\,\, \text{ is finite},
         \end{equation}
 and~\eqref{eq415} holds for all $n\leq m \leq 2n$, then the limit of $\{ a(n)/n\}$ still exists.

 Let us call a sequence  $\{ a(n)\}$  $(\mu, N, f)$-{\em subadditive} if~\eqref{eq415} holds for all $N\leq n\leq m \leq \mu n$.
 We usually suppose that $f$ is a non-negative monotone increasing real function but we will discuss more general cases as well.
Our Theorem~\ref{th32} yields the following corollary.

\begin{thm}\label{th44}
Suppose $\mu>1$ and $N\geq 1 $ are given and  $f$ is a non-negative monotone increasing real function.
If the sequence $\{ a(1)$, $a(2),\dots \}$ is $(\mu, N, f)$-subadditive, i.e.,
\[     a(n+m)\leq a(n)+a(m)+f(n+m) \quad \forall m\leq n\leq \mu m,\enskip m,n\geq N,
   \] 
then the $\lim_{n\to \infty}\dfrac{a(n)}{n}$ exists.
(It may be\enskip  $-\infty$).
\end{thm}

\smallskip
{\bf  Near subadditivity is {\em really} important}

Subadditivity is important, it appears in all parts of mathematics.
We all have our favorite examples and applications.
But nearly subadditivity is even more applicable,  here we mention a few areas.

In the beginning of the Bollob\'as-Riordan book~\cite{BR_2006}
 the  de Bruijn-Erd\H os theorem is listed  (as Lemma 2.1 on page 37) among the important useful tools
in Percolation Theory.
 The de Bruijn-Erd\H os theorem is widely used in investigating sparse random structures, e.g.,
 Bayati, Gamarnik, and Tetali~\cite{BGT_2009} (Proposition 5 on page 4011), Turova~\cite{T_2011}, 
or Kulczycki, 
Kwietniak, 
and Jian Li ~\cite{KKL} concerning  entropy of shift spaces.

Also, recurrence relations of type~\eqref{eq415} are often encountered in the analysis of divide and conquer algorithms,

\centerline{
${a(n+m)\leq a(n)+a(m) \enskip  +}$ {cost of cutting.} }

see, e.g.,  Hsien-Kuei Hwang and  Tsung-Hsi Tsai~\cite{HT_2003}.
In Economics it is an essential property of some {\em cost} functions that
 {COST(X+Y) $\leq$ COST(X)+COST(Y)}.
Similar relations appear in Physics and  in Combinatorial optimization  (see, e.g., Steele \cite{Steele}).

Also see, e.g.,  Capobianco~\cite{Ca}
concerning  cellular automatas, 
 Ceccherini-Silberstein, 
Coornaert, 
and F. Krieger\cite{CCK} 
for an analogue on cancellative {\em amenable semigroups}.

\smallskip
{\bf  Proof of Theorem~\ref{th44} using Theorem~\ref{th32}}

We utilize the proof from~\cite{deBE_1952} (bottom of page 163).
For $n\geq N$ define
\[
   G(n):= a(n) + 3n \left( \sum_{x\geq n}  f(x)/x^2 \right).
\]
Then the monotonicity of $f$, the relation $n\leq m\leq 2n$, and an easy calculation imply that
\[
   G(n+m) \leq G(n)+ G(m)
\]
whenever~\eqref{eq415} holds for $(n,m)$.

Theorem~\ref{th32} can be applied to $\{ G(n)\}$, so we have that the limit
\[
\lim_{n\to \infty} \frac{G(n)}{n}  = \lim_{n\to \infty}\left(  \frac{a(n)}{n} + \left( \sum_{x\geq n}  \frac{f(x)}{x^2} \right)\right)
\]
exists.
Here the last term tends to $0$ as $n\to \infty$ by~\eqref{eq416} and we are done.  \qed

\section{How large the error term $f(x)$ could be?}\label{S5}

It is very natural to ask how more one can extend the de Bruijn-Erd\H os theorem concerning $f$-nearly subadditive sequences
 (the case $\mu=2$, $N=1$).
Especially, how large the error term could be?

\smallskip
{\bf  $f(x)=o(x)$ is necessary}

Suppose that $f(n)$ is non-negative and $\limsup f(n)/n > L> 0$.
We can easily construct a sequence $\{ a(n)\}$ satisfying~\eqref{eq415} for all pairs $m,n\geq 1$
such that $\lim a(n)/n$ does not exist.
We do not even use that $f$ is monotone or not.

Given such an $f$ one can find a sequence of integers $1\leq n_1 < n_2 < n_3 < \dots$ such that
$f(n_{i})/n_i > L/2$,  and  $n_{i+1}\geq  n_i+2$ for all $i\geq 1$.
Define $a(n)= f(n_i)$ if $n=n_i$ and $0$ otherwise.  \qed

  \smallskip
{\bf  $f(x)=o(x)$ is not sufficient}

Condition~\eqref{eq416} allows $f(x)=O(x^{1-c})$ ($c> 0$ fixed)  or even $f(x)=O(x/(\log x)^{1+c})$.
The first author observed that $f(x)$ could not be $\Omega (x/\log x)$.
In 2016 he~\cite{FurediSch4} proposed the following problem for Schweitzer competition for university students (in Hungary).
``Prove that there exists a sequence  $a(1), a(2), \dots, a(n), \dots$ of real numbers such that
\begin{equation*}   a(n+m) \leq a(n)+ a(m) {+\frac{n+m}{\log(n+m)}}
   \end{equation*}
for all integers $m,n\geq 1$, and  the set
{$\{ a(n)/n: n\geq 1\}$ is everywhere dense} on the real line.''
(There were two correct solutions: by N\'ora Frankl, and Kada Williams and two partial solutions by Bal\'azs Maga,  and J\'anos Nagy).

\smallskip
{\bf deBruijn and Erd\H os got the best result}

  We show that the de Bruijn-Erd\H{o}s condition~\eqref{eq416} for the error term is not only sufficient but also necessary
 in the following strong sense.
\begin{thm}\label{th56}
Let   $f(n)$  be a non-negative, non-decreasing sequence and suppose
\begin{equation}\label{eq521} \sum_{1\leq n< \infty}  f(n)/n^2 = {\infty}.
\end{equation}
Then {there exists a nearly $f$-subadditive  sequence}  $b(1), b(2),  b(3), \dots$ of rational numbers, i.e., for all integers $m,n\geq 1$
\begin{equation*}   b(n+m) \leq b(n)+ b(m) { + f(n+m)}
   \end{equation*}
such that  the set of slopes  takes all rationals exactly once,
 ${ \{ b(n)/n: n\geq 1\}={\mathbf Q}}$.
\end{thm}

The proof is constructive and presented in the next section.

\section{Proof of Theorem~\ref{th56}, a construction}\label{S6}

A typical subadditive function is concave like, e.g.,  for $a(x)=\sqrt{x}$ we have $\sqrt{x+y}\leq \sqrt{x}+ \sqrt{y}$
 (for $x,y\geq 0$). The main idea of the construction for Theorem~\ref{th56} is that a nearly $f$-subadditive sequence $\{a(n)\}$ could be (strictly) convex with $\lim_{n\to \infty}   a(n)/n = \infty$.

\smallskip
{\bf A convex $f$-subadditive function}

\begin{claim}\label{cl61}
Suppose that $f(n)$ is a non-negative, non-decreasing sequence,
  $ 0\leq f(2)\leq f(3)\leq\dots$
Define $f(1)=a(1)=0$ and in general let
\begin{equation}\label{eq622}
       a(n):=   n \left( \sum_{i=1}^{n}   \frac{f(i)}{i^2}\right).
         \end{equation}
Then the sequence $\{ a(n) \}$ is nearly $f$-subadditive, it satisfies~\eqref{eq415}.
  \end{claim}

\noindent
{\sl Proof.}\enskip
Write down the definition of $a(n)$, simplify, use the monotonicity of $f$, finally the estimate
$\left( \sum_{u< i\leq v}   1/i^2 \right) < (1/u)-(1/v)$ (for integers  $1\leq u< v$).
 We obtain
\begin{eqnarray*}
&{}& a(n+m)-a(n)- a(m)\\
&{}&\quad\quad\quad\quad
=  n \left( \sum_{i\leq n+m}   \frac{f(i)}{i^2}\right)
  +m \left( \sum_{i\leq n+m}   \frac{f(i)}{i^2}\right)
  -n \left( \sum_{i\leq n}   \frac{f(i)}{i^2}\right)
  -m \left( \sum_{i\leq m}   \frac{f(i)}{i^2}\right)
\\
&{}&\quad\quad\quad\quad
=  n \left( \sum_{n< i\leq n+m}   \frac{f(i)}{i^2}\right)
  +m \left( \sum_{m<i\leq n+m}   \frac{f(i)}{i^2}\right)
\\
&{}&\quad\quad\quad\quad
 \leq  n  f(n+m)\left(   \frac{1}{n}- \frac{1}{n+m}\right)
  +m f(n+m)\left( \frac{1}{m}- \frac{1}{n+m}\right) = f(n+m).  \qed
\end{eqnarray*}

\begin{claim}\label{cl62}
The above sequence $\{ a(n)\}$ defined by~\eqref{eq622} is non-negative and convex,
 i.e., for $n\geq 2$ we have
\begin{equation*}
   a(n) \leq \frac{a(n-1)+ a(n+1)}{2}.
   \end{equation*}
  \end{claim}


\noindent
{\sl Proof.}\enskip
We have
\begin{eqnarray*}
a(n-1)+ a(n+1)-2 a(n) &=&  (n-1) \left( \sum_{i\leq n-1}   \frac{f(i)}{i^2}\right) +
 (n+1) \left( \sum_{i\leq n+1}   \frac{f(i)}{i^2}\right)-2 n \left( \sum_{i\leq n}   \frac{f(i)}{i^2}\right)  \notag \\
&= & \frac{f(n+1)}{(n+1)}-  (n-1)  \frac{f(n)}{n^2}
     \geq \frac{f(n+1)}{(n+1)n^2}\geq 0. \quad  
     \qed
\end{eqnarray*}

\smallskip
{\bf The end of the proof of Theorem~\ref{th56}}

In this section $\{ f(n)\}$ is given by Theorem~\ref{th56}, and $\{ a(n)\}$ is the well-defined nearly $f$-subadditive, convex sequence obtained by~\eqref{eq622} in Claim~\ref{cl62}. Then~\eqref{eq521} implies  $\lim_{n\to \infty}   {a(n)/n} = \infty$.

For the rest of the proof  the main observation is the following: If $c(1) \leq c(2)\leq c(3)\leq \dots$\enskip  is a monotone sequence, and $\{ a(n)\}$ is $f$-subadditive, then 

\medskip
\centerline{$b(n):= a(n)-c(n)n$ \enskip is $f$ subadditive as well.}

${}$
\vskip -7mm
 Indeed,
${}$
\vskip -9mm
\begin{multline*}b(n+m)-b(n)-b(m)-f(n+m) \\= [a(n+m)- c(n+m)(n+m)]-[a(n)-c(n)n]-[a(m)-c(m)m]-f(n+m)
\\ = [a(n+m)-a(n)-a(m)-f(n+m)] +(c(n)- c(n+m))n +(c(m)- c(n+m))m \leq 0.
  \end{multline*}

Let $r_1, r_2, r_3, \dots$ be an enumeration of $\mathbf Q$.
We will define a sequence  $1\leq n_0\leq n_1\leq n_2\leq \dots$ and simultaneously $\{ c(n)\}$  (and thus $\{ b(n)\}$ as well) such that
\newline
${}$ (D) 
\quad the slopes $\{ b(n)/n: 1\leq n\leq n_i\}$ are all distinct and rational, and 
\newline
${}$${}$ (R) 
\quad $r_i\in \{ b(n)/n: 1\leq n\leq n_i\}$,  \enskip($i\geq 1$). 

We proceed by induction on $i$.
Let $n_0$ be the smallest $x\geq 1$ such that $f(x)>0$. Equation~\eqref{eq521} implies that $1\leq n_0< \infty$.
Choose $c(1)\leq  \dots\leq  c(n_0)$ arbitrarily such that the fractions
 $b(x)/x =$ $(a(x)-c(x)x)/x$ are all rationals and they are all distinct.
Since these are finitely many constraints of the form
\[  \dfrac{a(x)}{x}-c(x) \neq  \dfrac{a(y)}{y}-c(y)\quad 1\leq x\neq y\leq n_0
\]
and the set $\mathbf Q$ is everywhere dense on $\mathbf R$, one can easily choose appropriate $c(x)$'s.

If $n_0, n_1, \dots, n_i$ has been already defined (satisfying properties (D) and (R)) then proceed as follows.

If $r_{i+1}\in \{ b(x)/x: 1\leq x \leq n_i\}$, then let $n_{i+1}:= n_i$.

If $r_{i+1}\notin  \{ b(x)/x: 1\leq x \leq n_i\}$ then define $n_{i+1}$ as the smallest integer $x$
 satisfying
\[   x> n_i, \quad \dfrac{a(x)}{x}-c(n_{i})  > r_{i+1}.
\]
Such $x$ exists. Let $c(n_{i+1}):= \dfrac{a(n_{i+1})}{n_{i+1}}- r_{i+1}$.
It follows that $c(n_i)< c(n_{i+1})$.
Then select $c(x)$ for integers $x$ with $n_i< x < n_{i+1}$ such that
 the values of $a(x)/x  -c(x) $ are all rationals, distinct from each other,
 have no common values with $\{ b(n)/n: 1\leq n\leq n_i\}\cup \{ r_{i+1}\}$ and also
$c(n_i) \leq c(n_i+1 ) \leq \dots \leq c(n_{i+1})$.
These are finitely many conditions but $c(n_i)< c(n_{i+1})$ and $\mathbf Q$ is everywhere dense, so the induction step can be done. This completes the construction. \qed

\section{Conclusion, problems}\label{S7}

Let $X\subseteq {\mathbf N}\times {\mathbf N}$, $f: {\mathbf N}\to {\mathbf R}$.
The sequence $\{ a(n)\}$ is $(X,f)$-subadditive if $a(m+n)\leq a(n)+a(m)+f(n+m)$ holds for
 $(n,m)\in X$.
We have found conditions for $X$ and $f$, strengthening the original Fekete's lemma and its de Bruijn-Erd\H os generalization,  which ensure that $\lim a(n)/n$ exists.
Certainly further thinning of $X$ are possible. We mention two of these problems. 

Is it possible to replace the constraint $n\leq m \leq \mu n$ in Theorem~\ref{th32} by some condition like  $n\leq m \leq n+r(n)$ where $r(n)=o(n)$ some slow growing function? (Probably not).

What is the structure of $1^+$ subadditive sequences? Can we tell more than Lemma~\ref{le33}?

Finally, it is well-known that if $a(x)$ is a {\em measurable} subadditive function $a:(0, \infty)\to {\mathbf R}$, then the limit $\lim_{x\to \infty} a(x)/x $ exists.
The non-measurable subadditive functions include the {\em Cauchy-functions} which do not have limits, and are far from linear.
This is a large field of analysis, and alos in number theory concerning  additive functions.
There are many results and questions, see, e.g.,~\cite{Er46, Ru, Wi}.

\end{document}